\documentclass{article}
\usepackage{amsmath}
\usepackage{ntheorem-hyper}
\usepackage{amssymb}
\usepackage{authblk} % author envir package
\usepackage{algorithm}
\usepackage{algorithmicx}
\usepackage{algpseudocode}
\usepackage{graphicx}
\usepackage{graphbox}
\usepackage{subfigure}
\usepackage{svg}
\usepackage{epstopdf}
\usepackage{tabularx}  
\usepackage{booktabs}
\usepackage[numbers, sort]{natbib}
\usepackage{hyperref}
\usepackage[top=2cm, bottom=2cm, left=2cm, right=2cm]{geometry}  
\usepackage{xcolor} % \textcolor{red/blue/green/black/white/cyan/magenta/yellow/teal}{text}
\usepackage{color}
\usepackage[title]{appendix}

\usepackage[T1,T2A]{fontenc}
\usepackage[utf8]{inputenc}
\usepackage[russian, german, english]{babel}

\hypersetup{urlcolor=blue, citecolor=red}

\newtheorem{theorem}{Theorem}
\newtheorem{prop}{Proposition}
\theorembodyfont{\upshape}
\newtheorem*{proof}{Proof}

\newtheorem{corollary}{Corollary}
\newtheorem{lemma}{Lemma}

\def\RR{\mathbb{R}}

\begin{document}
\title{Riemannian convex bodies with a boundary geodesic segment}
\author{Hongda Qiu}
\date{}
\maketitle

\begin{abstract}
We construct compact convex bodies whose boundaries contain a nontrivial geodesic segment and are strictly convex elsewhere. The construction applies if sectional curvature is constant near a point, or if a normal Jacobi operator has a simple maximal eigenvalue. This implies existence of such bodies in any three-dimensional manifold. A further argument to show existence in dimension four.
%The proof combines a
%short-interval criterion for Jacobi convexity of parallel cones, a
%Euclidean body with arbitrarily skinny transverse tangent cones, and
%uniform Hessian estimates in coordinates defined by a parallel normal
%frame and the normal exponential map.
\end{abstract}

\section{Introduction}

A subset of a Riemannian manifold
is called \emph{convex} if every minimizing geodesic between any two of
its points is contained in this subset. Alexander~Lytchak and Anton~Petrunin~\cite[Question~B.1]{LP} asked whether every Riemannian manifold contains a nontrivial geodesic lying in the
boundary of a convex subset. We establish the following sufficient
condition, which gives an affirmative answer in dimension three.
All manifolds are smooth and connected.
%completeness is not assumed.

\begin{theorem}
\label{thm:main}
Let $(M^n,g)$, $n\geq3$, be a Riemannian manifold, let $p\in M$, and let
$v\in T_pM$ be a unit vector and $v^{\perp}$ its normal space. Assume that either
\begin{enumerate}
\item $g$ has constant sectional curvature in a neighborhood of $p$; OR
\item the normal Jacobi operator
\[
A_v:v^\perp\longrightarrow v^\perp,\qquad w\mapsto R(w,v)v,
\]
has a simple maximal eigenvalue.
\end{enumerate}
Then every neighborhood of $p$ contains a compact convex body
$\mathcal C$ with nonempty interior such that
\[
\gamma([-\rho,\rho])\subset\partial\mathcal C,
\qquad \gamma(t)=\exp_p(tv),\quad \rho>0.
\]
Specifically, the boundary is smooth and has a positive-definite second fundamental
form outside this segment; every pair of points of $\mathcal C$ is
joined by a unique minimizing geodesic in $M$, and this geodesic is
contained in $\mathcal C$;
the transverse tangent cones along the relative
interior of the segment may be chosen arbitrarily skinny.
\end{theorem}

Here, we say that a transverse cone is \emph{$\epsilon$-skinny} if its nonzero
vectors make an angle at most $\arctan\epsilon$ with a fixed unit vector. For Theorem 1, we use the following cones in an orthonormal frame:
\begin{equation}
\label{eq:skinny-cone}
\mathop{\mathrm{Cone}}\nolimits_\epsilon
=\left\{\sum_{i=2}^n x_iE_i:
x_2\geq0,\quad \sqrt{x_3^2+\cdots+x_n^2}\leq\epsilon x_2\right\}.
\end{equation}
Their parallel translates, together with \(\mathbb R\dot\gamma(t)\), are the tangent cones along the relative interior of the segment.
For every fixed $\epsilon>0$, the construction allows arbitrarily small
$\rho$; no uniform lower bound for $\rho$ as $\epsilon\downarrow0$ is asserted.

\begin{corollary}
\label{cor:three}
Every three-dimensional Riemannian manifold contains a compact
convex body with nonempty interior whose boundary contains a nontrivial
geodesic segment and is smooth and strictly convex elsewhere.
\end{corollary}

\begin{proof}
In dimension three, a nonscalar normal Jacobi operator has a simple maximal eigenvalue. If all normal Jacobi operators are scalar, sectional curvature is pointwise constant and hence constant by Schur’s theorem. Theorem~\ref{thm:main} applies in either case.
\end{proof}

We also settle the four-dimensional case of Lytchak and Petrunin's problem. In this construction, the boundary may be nonsmooth outside the geodesic segment.
\begin{theorem}
	\label{thm:four}
	Every smooth four-dimensional Riemannian manifold contains a compact
	convex body with nonempty interior whose boundary contains a nontrivial
	geodesic segment. Moreover, the transverse tangent cones along the relative
	interior of this segment can be chosen to be parallel, pointed, and
	arbitrarily skinny.
\end{theorem}

The following case remains open in higher dimensions.
%[Repeated maximal Jacobi eigenvalue]

\textbf{Open Problem.} Let $M$ be a connected smooth Riemannian manifold of dimension
$n\geq5$ with nonconstant sectional curvature. Suppose that, for
every $p\in M$ and every unit vector $v\in T_pM$, the maximal
eigenvalue of the normal Jacobi operator $ A_v:v^\perp\longrightarrow v^\perp,
\quad A_v(w)=R(w,v)v $ has multiplicity at least two. Does there exist a compact convex body $\mathcal{C}\subset M$ with nonempty
interior and a nonconstant geodesic segment
$\gamma:[a,b]\to M$ such that $\gamma([a,b])\subset\partial \mathcal{C}$?

For the rest of content, we use the curvature convention for which the Jacobi equation is
$D_t^2J+R(J,\dot\gamma)\dot\gamma=0$. For an inward unit normal
$\nu$, the convex sign convention is
$\mathrm{II}(X,Y)=g(-\nabla_X\nu,Y)$.

\section{Jacobi convexity of parallel skinny cones}

Let $\gamma:(-l,l)\to M$ be a unit-speed geodesic and choose a parallel
orthonormal frame $E_1=\dot\gamma,E_2,\ldots,E_n$ along it. We identify every normal space \(\dot\gamma(t)^\perp\) with the fixed space \(E=\dot\gamma(0)^\perp\) by parallel transport along \(\gamma\).
Write $A(t)$ for the normal Jacobi operator in this identification;
normal Jacobi fields then satisfy $\ddot{J}+A(t)J=0$.

%A parallel family of closed convex cones $K_t\subset\dot\gamma(t)^\perp$
%is called \emph{Jacobi convex} if every normal Jacobi field whose
%endpoint values on a subinterval $[a,b]$ lie in $K_a,K_b$ lies in
%$K_t$ throughout that subinterval. The intervals considered below
%are short enough that the Jacobi boundary-value problem is uniquely
%solvable. Equivalently, this is the corresponding property of the
%union of the cones in the normal bundle. The same property holds
%for $\RR\dot\gamma(t)\oplus K_t$, since the tangential component of
%a Jacobi field is affine.
We say that a parallel family of closed convex cones $K_t\subset\dot\gamma(t)^\perp$ is \textit{Jacobi convex} if for every subinterval $[a,b]\subset(-l,l)$, every normal Jacobi field with endpoint values in $K_a,K_b$ remains in $K_t$ throughout $[a,b]$. Note that if $K_t$ is Jacobi convex, then so is $\mathbb R\dot\gamma(t)\oplus K_t$.

%We consider short enough intervals on which the Jacobi boundary-value problem is uniquely solvable.

\begin{lemma}
\label{lem:jacobi}
Let $K\subset E$ be a closed convex cone, with inward dual
$K^*=\{\nu:\langle V,\nu\rangle\geq0\text{ for every }V\in K\}$.
Suppose that
\begin{equation}
\label{eq:inward}
\langle A(t)V,\nu\rangle\geq0
\quad\text{if }V\in K,\ \nu\in K^*,\ \langle V,\nu\rangle=0.
\end{equation}
If $[a,b]\subset(-l,l)$ and
\begin{equation}
\label{eq:short}
\left(\sup_{[a,b]}\|A(t)\|_{\mathrm{op}}\right)(b-a)^2<\pi^2,
\end{equation}
then the unique Jacobi field with any prescribed endpoint values in
$K$ remains in $K$ on $[a,b]$.
\end{lemma}

\begin{proof}
Suppose  $ \ddot{Y}+AY=0,\quad Y(a)=Y(b)=0$. Integration by parts and the Dirichlet Poincar\'e
inequality~\cite[Theorem~257, p.~185]{HLP} give
\[
\int_a^b|\dot{Y}|^2
=\int_a^b\langle AY,Y\rangle
\leq \sup_{[a,b]}\|A(t)\|_{\mathrm{op}}\cdot\int_a^b|Y|^2
\leq\frac{\sup_{[a,b]}\|A(t)\|_{\mathrm{op}}\cdot(b-a)^2}{\pi^2}\int_a^b|\dot{Y}|^2.
\]
Here the scalar inequality on $[0,\pi]$ is rescaled to $[a,b]$
and applied componentwise.
Thus $Y=0$. For solutions of the Jacobi equation with \(Y(a)=0\), the map sending \(\dot Y(a)\) to \(Y(b)\) is linear and injective, hence bijective. This gives existence and uniqueness for arbitrary endpoint values by superposition.

Let $J$ be any solution to $ \ddot J+A(t)J=0$ with $ J(a),J(b)\in K $. Write $J_{\parallel}=\operatorname{proj}_K J$ and $J_{\perp}=J-J_{\parallel}$.
The projection inequality
$\langle J_{\perp},V-J_{\parallel}\rangle\leq0$ for any $V\in K$ implies
$\langle J_{\perp},J_{\parallel}\rangle=0$ and $-J_{\perp}\in K^*$. Consequently,
\begin{equation}
\label{eq:projection-sign}
\langle AJ_{\parallel},J_{\perp}\rangle\leq0.
\end{equation}
Applying the projection inequality at two parameter values gives
$\langle J_{\parallel}(t)-J_{\parallel}(s),J_{\perp}(t)-J_{\perp}(s)\rangle\geq0$.
The projection is Lipschitz, so
$\langle \dot{J}_{\parallel},\dot{J}_{\perp}\rangle\geq0$ almost everywhere and hence
$\langle \dot{J},\dot{J}_{\perp}\rangle\geq|\dot{J}_{\perp}|^2$.
Since $J_{\perp}(a)=J_{\perp}(b)=0$, we obtain
\[
\begin{aligned}
\int_a^b|\dot{J}_{\perp}|^2
&\leq\int_a^b\langle \dot{J},\dot{J}_{\perp}\rangle
=\int_a^b\langle AJ,J_{\perp}\rangle\\
&\leq \sup_{[a,b]}\|A(t)\|_{\mathrm{op}}\cdot\int_a^b|J_{\perp}|^2
\leq\frac{\sup_{[a,b]}\|A(t)\|_{\mathrm{op}}\cdot (b-a)^2}{\pi^2}\int_a^b|\dot{J}_{\perp}|^2.
\end{aligned}
\]
Condition~\eqref{eq:short} forces $J_{\perp}=0$, as required.

\end{proof}

Condition~\eqref{eq:inward} expresses the inward curvature condition
of~\cite[Key Lemma~2.1(i)]{LP} in terms of supporting normals.
Lemma~\ref{lem:jacobi} provides the short-interval
sufficiency statement for a prescribed parallel family of cones.
If $M$ has constant sectional curvature $\lambda$, then $A(t)=\lambda
\operatorname{Id}$ and the left side of~\eqref{eq:inward} is zero.
Thus every parallel family of closed convex cones is Jacobi convex
on a sufficiently short geodesic, independently of the aperture
of the cones. We next consider the simple maximal eigenvalue case.

\begin{prop}
\label{prop:jacobi}
Suppose $A(0)$ has eigenvalues $\lambda_2,\ldots,\lambda_n$ with
$\lambda_2>\max_{i\geq3}\lambda_i$, and choose $E_2(0),\ldots,E_n(0)$
as the corresponding orthonormal eigenbasis. For every $\epsilon>0$,
the parallel copies of $\mathop{\mathrm{Cone}}\nolimits_\epsilon$ in~\eqref{eq:skinny-cone} are
Jacobi convex after sufficiently shortening $\gamma$. 
\end{prop}
%Specifically, this cone is pointed, has
%nonempty interior, and has smooth boundary away from its vertex.

\begin{proof}
Put $\delta=\lambda_2-\max_{i\geq3}\lambda_i>0$. For a unit vector
$u\in\operatorname{span}\{E_3(0),\ldots,E_n(0)\}$, set
\[
q_u=E_2(0)+\epsilon u,\qquad \nu_u=\epsilon E_2(0)-u.
\]
The cone is the intersection of the half-spaces
$\langle V,\nu_u\rangle\geq0$, and its contact face with
$\nu_u^\perp$ is $\RR_{\geq0}q_u$. At $t=0$,
\[
\langle A(0)q_u,\nu_u\rangle
=\epsilon\bigl(\lambda_2-\langle A(0)u,u\rangle\bigr)
\geq\epsilon\delta.
\]
Since $|q_u||\nu_u|=1+\epsilon^2$, continuity gives, uniformly in $u$,
\[
\langle A(t)q_u,\nu_u\rangle\geq\epsilon\delta/2>0
\]
on a sufficiently short interval. Every nonzero supporting normal
at a nonzero boundary point is a positive multiple of the
corresponding $\nu_u$. At the vertex, Condition~\eqref{eq:inward} is automatic;
at an interior point, the only orthogonal dual vector is zero.
Thus Lemma~\ref{lem:jacobi} applies after shortening the interval
again to ensure~\eqref{eq:short}. 

%The remaining assertions follow
%directly from~\eqref{eq:skinny-cone}: the cone is an epigraph of a
%norm, its intersection with its negative is $\{0\}$, and its nonzero boundary is the regular level set
%$\epsilon^2x_2^2-\sum_{i\geq3}x_i^2=0$.

\end{proof}

\section{A convex body with skinny transverse tangent cones in $\RR^n$}

Let $e_1,\ldots,e_n$ denote the standard coordinate frame of $\RR^n$,
identified with the parallel frame used above. The \textit{tangent cone} of a convex body at a boundary point is the local Hausdorff limit obtained by translating that point to the origin and dilating the body by factors converging to infinity.
\begin{prop}
		\label{prop:closed-convex-hypersurface-skinny-cones}
		Let $n\geq3$ and $\epsilon>0$. There exists a compact convex body
		$C_\epsilon\subset\RR^n$ with nonempty interior whose boundary
		$\Sigma_\epsilon=\partial C_\epsilon$ contains the segment
		\[
		I=\{(x_1,0,\ldots,0):-1\leq x_1\leq1\}.
		\]
		Moreover, 
		\begin{enumerate}
			\item The hypersurface $\Sigma_\epsilon$ is smooth outside $I$. 
			\item At every point
			$p\in I\setminus\{(-1,0,\ldots,0),(1,0,\ldots,0)\}$, the tangent cone,
			translated to the origin, is
			\begin{equation}
				\label{eq:tangent-cones}
				\begin{aligned}
					T_pC_\epsilon
					&=\{(x_1,\ldots,x_n)\in\RR^n:
					x_2\geq0,\ \sqrt{x_3^2+\cdots+x_n^2}\leq\epsilon x_2\}.
				\end{aligned}
			\end{equation}
%			,\\
%			T_p\Sigma_\epsilon
%			&=\{(x_1,\ldots,x_n)\in\RR^n:
%			x_2\geq0,\ \sqrt{x_3^2+\cdots+x_n^2}=\epsilon x_2\}
			In particular, its transverse section \(T_pC_\epsilon\cap\{x_1=0\}\) is pointed, convex, and \(\epsilon\)-skinny.
			\item $\Sigma_\epsilon$ has a positive-definite second fundamental form at
			every point of $\Sigma_\epsilon\setminus I$.
		\end{enumerate}
	\end{prop}
	
	\begin{proof}
		Consider a closed plane curve $f(x,y)=0$ (Figure~\ref{fig:graph of f}), where
		\[
		f(x,y)=x^2-\epsilon^2(y^2-y^3),\qquad y\geq0.
		\]
		\begin{figure}[htbp]
			\centering
			\includegraphics[width=0.25\textwidth]{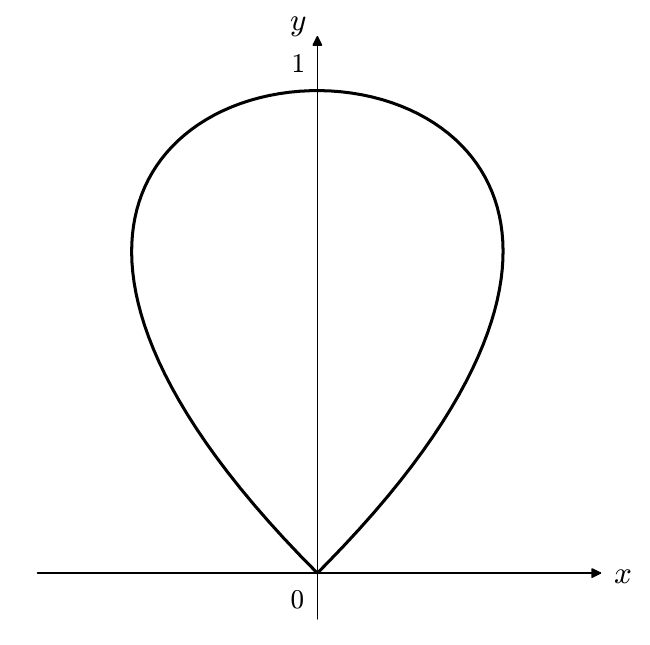}
			\caption{$f(x,y)=0, y\geq 0$.}
			\label{fig:graph of f}
		\end{figure}
		
		Define the hypersurface directly by
		\begin{equation}
			\label{eq:prescribed-hypersurface}
			\Sigma_\epsilon
			=\overline{\left\{(x_1,\ldots,x_n)\in\RR^n:
				\begin{array}{l}
					-1<x_1<1,\quad x_2\geq0,\\[2pt]
					\displaystyle
					f\left(\frac{\sqrt{x_3^2+\cdots+x_n^2}}{1-x_1^2},
					\frac{x_2}{1-x_1^2}\right)=0
				\end{array}\right\}},
		\end{equation}
		where the closure is taken in $\RR^n$. Note that the closure adds precisely the two	points $(-1,0,\ldots,0)$ and $(1,0,\ldots,0)$. Specify its enclosed region by
		\begin{equation}
			\label{eq:convex-body}
			C_\epsilon
			=\overline{\left\{(x_1,\ldots,x_n)\in\RR^n:
				\begin{array}{l}
					-1<x_1<1,\quad x_2\geq0,\\[2pt]
					\displaystyle
					f\left(\frac{\sqrt{x_3^2+\cdots+x_n^2}}{1-x_1^2},
					\frac{x_2}{1-x_1^2}\right)\leq0
				\end{array}\right\}}.
		\end{equation}
		For $|x_1|<1$, the equation for $\Sigma_\epsilon$ is equivalent to
		\begin{equation}
			\label{eq:hypersurface-in-slab}
			x_3^2+\cdots+x_n^2
			=\epsilon^2x_2^2\left(1-\frac{x_2}{1-x_1^2}\right),
			\qquad 0\leq x_2\leq1-x_1^2.
		\end{equation}

		\begin{figure}[htbp]
			\centering
			\includegraphics[width=0.5\textwidth]{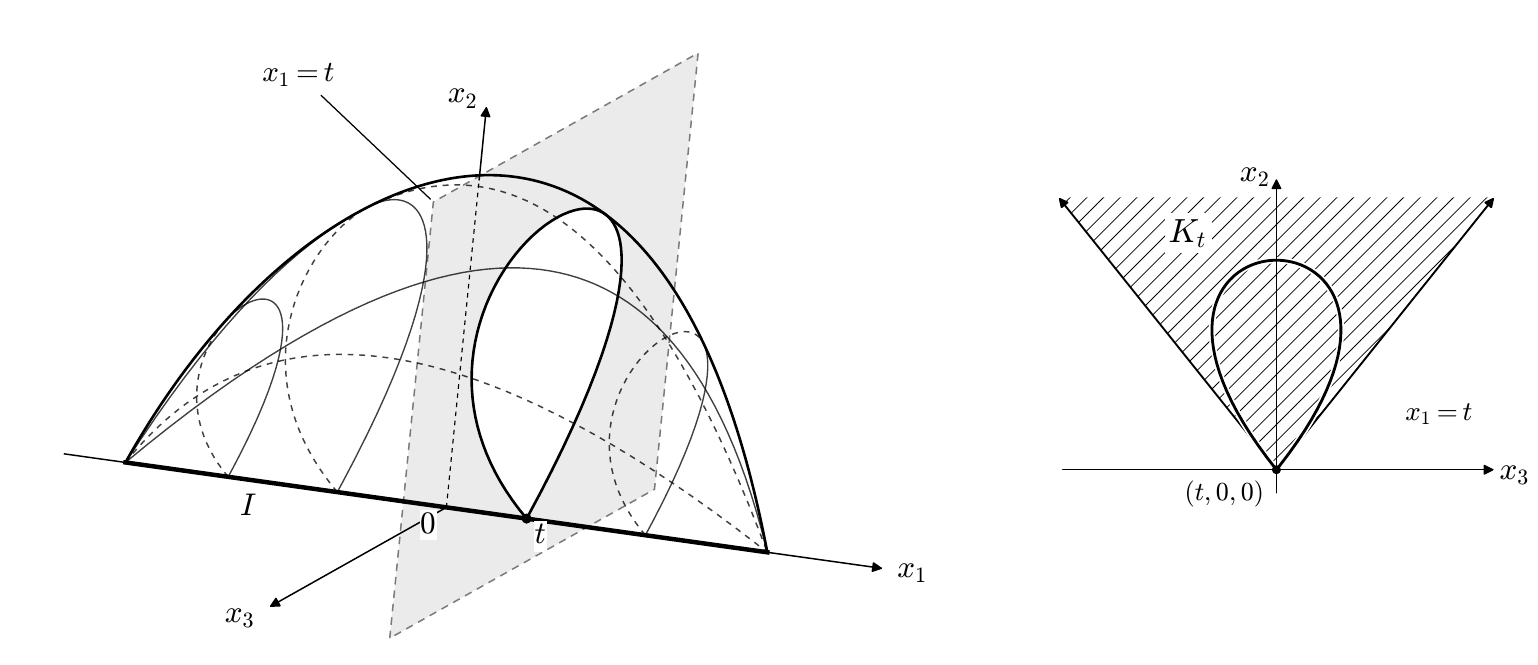}
			\caption{An illustration in dimension three: the body $C_\epsilon$ and the normal plane $x_1=t$ (left); its section and transverse tangent cone $K_t$ at $(t,0,0)$ (right).}
			\label{fig:graph of tangent cone}
		\end{figure}
		
		We prove convexity directly from this equation. On the convex domain
		\[
		D=\{(x_1,x_2)\in\RR^2:x_2\geq0,\ x_1^2+x_2\leq1\},
		\]
		define the radius function
		\[
		G(x_1,x_2)
		=\epsilon x_2\sqrt{1-\frac{x_2}{1-x_1^2}}
		\qquad(|x_1|<1),
		\]
		and set $G(-1,0)=G(1,0)=0$. This extension is continuous because
		$0\leq G(x_1,x_2)\leq\epsilon x_2$. In the interior of $D$,
		a direct calculation gives
		\[
		\frac{\partial^2G}{\partial x_2^2}
		=\frac{\epsilon\bigl(3x_2-4(1-x_1^2)\bigr)}
		{4(1-x_1^2)^2\left(1-\dfrac{x_2}{1-x_1^2}\right)^{3/2}}<0
		\]
		and
		\[
		\det D^2G
		=\frac{\epsilon^2x_2^2\bigl(4(1-x_1^2)-3x_2\bigr)}
		{4(1-x_1^2)^2(1-x_1^2-x_2)^2}>0.
		\]
		Thus $D^2G$ is negative definite in the interior of $D$, and continuity
		implies that $G$ is concave on all of $D$. Since
		\begin{equation}
			\label{eq:body-in-slab}
			C_\epsilon
			=\{(x_1,\ldots,x_n):(x_1,x_2)\in D,\
			\sqrt{x_3^2+\cdots+x_n^2}\leq G(x_1,x_2)\},
		\end{equation}
		the triangle inequality and concavity of $G$ imply that $C_\epsilon$ is
		convex. It is compact, and $(0,\tfrac12,0,\ldots,0)$ is an interior point.	Furthermore,
		\[
		C_\epsilon\subseteq\{x_2\geq0\},
		\qquad C_\epsilon\cap\{x_2=0\}=I,
		\]
		so $I\subset\Sigma_\epsilon$.
		
		We now compute the tangent cones of \(C_\epsilon\) along the relative interior of \(I\). Fix $p=(x_1^0,0,\ldots,0)$ with $-1<x_1^0<1$. Under the rescaling
		\[
		(x_1,\ldots,x_n)
		\longmapsto(x_1^0+h x_1,h x_2,\ldots,h x_n),\qquad h>0,
		\]
		the inequality defining $C_\epsilon$ becomes
		\[
		x_2\geq0,\qquad
		\sqrt{x_3^2+\cdots+x_n^2}
		\leq\epsilon x_2
		\sqrt{1-\frac{h x_2}{1-(x_1^0+h x_1)^2}},
		\]
		with the additional domain restrictions automatic for each fixed direction
		and all sufficiently small $h>0$. The right-hand side tends to
		$\epsilon x_2$ as $h\downarrow0$. Every direction satisfying the strict
		inequality $\sqrt{x_3^2+\cdots+x_n^2}<\epsilon x_2$ is therefore realized
		for all sufficiently small $h$. Conversely, the original inequality
		always implies $\sqrt{x_3^2+\cdots+x_n^2}\leq\epsilon x_2$.
		Taking closures proves the identity in~\eqref{eq:tangent-cones}.
		
%		Rescaling~\eqref{eq:hypersurface-in-slab} shows that every vector in
%		$T_p\Sigma_\epsilon$ satisfies $x_2\geq0$ and
%		$\sqrt{x_3^2+\cdots+x_n^2}=\epsilon x_2$. Conversely, for any such vector
%		with $x_2>0$, the points
%		\[
%		p+h x_1e_1+h x_2e_2
%		+h\sqrt{1-\frac{h x_2}{1-(x_1^0+h x_1)^2}}
%		\sum_{i=3}^n x_i e_i
%		\]
%		lie on $\Sigma_\epsilon$ for all sufficiently small $h>0$ and have the
%		prescribed tangent direction at $p$. When $x_2=0$, the remaining transverse
%		coordinates vanish, and the direction is realized by points of $I$.
%		This proves the second identity in~\eqref{eq:tangent-cones}.
		
		Finally, the hypersurface outside $I$ is given by \(F=0\) on \(|x_1|<1,\ x_2>0\), where
		\[
		F(x_1,\ldots,x_n)
		=x_3^2+\cdots+x_n^2-\epsilon^2x_2^2
		+\frac{\epsilon^2x_2^3}{1-x_1^2}.
		\]
		If some $x_i\neq0$ for $i\geq3$, then $\partial F/\partial x_i=2x_i\neq0$.
		Otherwise, outside $I$ the boundary equation gives $x_2=1-x_1^2>0$, and
		\[
		\frac{\partial F}{\partial x_2}=\epsilon^2(1-x_1^2)>0.
		\]
		The implicit function theorem proves that $\Sigma_\epsilon$ is smooth
		outside $I$.
		
		To prove positive definiteness of the second fundamental form outside $I$, use the convention
		$\mathrm{II}(V,V)=D^2H[V,V]/|\nabla H|$ for an outward defining function
		$H$. Where $x_3^2+\cdots+x_n^2>0$, take
		\[
		H(x_1,\ldots,x_n)=\sqrt{x_3^2+\cdots+x_n^2}-G(x_1,x_2).
		\]
		Let $V$ be a nonzero tangent vector, with components $V_1,\ldots,V_n$.
		The negative definiteness of $D^2G$ and the positive semidefiniteness of
		the Hessian of the Euclidean norm give $D^2H[V,V]>0$ whenever
		$(V_1,V_2)\neq(0,0)$. If $V_1=V_2=0$, tangency gives
		$\sum_{i=3}^n x_iV_i=0$, and hence
		\[
		D^2H[V,V]
		=\frac{\sum_{i=3}^n V_i^2}{\sqrt{x_3^2+\cdots+x_n^2}}>0.
		\]
		At the remaining points outside $I$, we have $x_2=1-x_1^2>0$ and
		$x_3=\cdots=x_n=0$. Tangency to $F=0$ gives $V_2=-2x_1V_1$, and direct
		differentiation yields
		\[
		D^2F[V,V]
		=2\epsilon^2(1-x_1^2)V_1^2+2\sum_{i=3}^n V_i^2>0
		\]
		for every nonzero tangent vector $V$. Since $F$ is an outward defining
		function with nonvanishing gradient, the second fundamental form is
		positive definite at these points as well.
		
	\end{proof}

\section{Local convexity}

Let $\gamma(0)=p$ and choose the parallel frame as in Section~2,
with an eigenbasis at $p$ in the simple-eigenvalue case. On a
neighborhood $U\subset\RR^n$ of the axis segment, consider
\begin{equation}
\label{eq:fermi}
\Phi(t,y_2,\ldots,y_n)
=\exp_{\gamma(t)}\left(\sum_{i=2}^n y_iE_i(t)\right).
\end{equation}
%After restricting $U$, this is a diffeomorphism onto a tubular
%neighborhood of $\gamma$. These are Fermi coordinates: the frame is
%parallel along $\gamma$, and the transverse coordinates are supplied
%by the normal exponential map. In particular, the manifold body
%associated with $\rho C_\epsilon\subset U$ is
%$\Phi(\rho C_\epsilon)$.
After shrinking $U$, the map $\Phi$ is a diffeomorphism onto a tubular neighborhood of $\gamma$. We use the pullback metric $\Phi^*g$ on $U$.

\begin{prop}
\label{prop:local}
Under either hypothesis of Theorem~\ref{thm:main}, let $C_\epsilon$
be the body of Proposition~\ref{prop:closed-convex-hypersurface-skinny-cones}. For each fixed
$0<\epsilon\leq1$, there is $\rho_1>0$ such that, for
$0<\rho<\rho_1$, the body
\[
\mathcal C_{\epsilon,\rho}=\Phi(\rho C_\epsilon)
\]
is locally geodesically convex at every point, including the endpoints
of $\gamma([-\rho,\rho])$. Specifically, its boundary is smooth and has a positive-definite
second fundamental form outside this segment. Here local convexity
means that each point has a strongly convex neighborhood $B$ for
which $B\cap\mathcal C_{\epsilon,\rho}$ is geodesically convex.
\end{prop}

\begin{proof}
Fix $0<\epsilon\leq1$ and take the body constructed in
Proposition~\ref{prop:closed-convex-hypersurface-skinny-cones}. Recall the function $F$ from the preceding proof. 
Write $\mathbf{x} =(x_1,\ldots,x_n) = (t, \mathbf{y})$ and
$\theta=x_2/(1-t^2)$. We sometimes regard normal vectors as vectors in $\RR^n$ with zero first coordinate.

\emph{Estimates of $F$.} We first estimate the Hessian of $F$ on tangent vectors to the regular boundary. On $|t|<1$, $x_2>0$, we have
\begin{equation}
\label{eq:model-function}
F(t,\mathbf{y})=\sum_{i=3}^n x_i^2-\epsilon^2x_2^2
+\frac{\epsilon^2x_2^3}{1-t^2}.
\end{equation}
The body is locally $\{F\leq0\}$, and $dF\ne0$ on
$\Sigma_\epsilon\setminus I$, as proved above. Write \(\partial_i=\partial/\partial x_i\). For a tangent vector $\mathbf{X}=\sum_{i=1}^n X_i\partial_i$, let
$\mathbf{X}_{\perp}=\sum_{i=2}^n X_i\partial_i$ be its transverse component. All coordinate norms below are Euclidean. Recalling $F=0$ and $dF(\mathbf{X})=0$, we have
\begin{equation}
\label{eq:squares}
\begin{aligned}
D^2F(\mathbf{X},\mathbf{X})
&=2\sum_{i=3}^n\left(X_i-\frac{x_i}{x_2}X_2\right)^2\\
&\quad+2\epsilon^2\theta(X_2+2t\theta X_1)^2
+2\epsilon^2x_2\theta^2X_1^2.
\end{aligned}
\end{equation}
%In particular, this identity is asserted for boundary tangent vectors.
On the regular boundary one also has
\begin{equation}
\label{eq:model-identities}
\begin{gathered}
\sum_{i=3}^n x_i^2=\epsilon^2x_2^2(1-\theta),\qquad
d_\mathbf{y}F(\mathbf{y})=\epsilon^2x_2^2\theta,\qquad \partial_i{F}=2x_i\ (i\geq3),\\
|dF|\leq4\epsilon x_2,\qquad |\mathbf{y}|\leq\sqrt2x_2,\qquad
0<x_2\leq\theta\leq1.
\end{gathered}
\end{equation}
The gradient bound follows directly from
$\partial_1{F}=2\epsilon^2t x_2\theta^2$,
$\partial_2{F}=\epsilon^2x_2(3\theta-2)$, and $\partial_i{F}=2x_i$ for $i\geq3$.
Moreover,
\begin{equation}
\label{eq:lower-bound}
D^2F(\mathbf{X},\mathbf{X})\geq\frac{2\epsilon^2x_2}{25}
\bigl(|\mathbf{X}_{\perp}|^2+\theta^2X_1^2\bigr).
\end{equation}
To see this, let $P=\sum_{i\geq3}(X_i-x_i X_2/x_2)^2$,
$Q=(X_2+2t\theta X_1)^2$, and $u=\theta X_1$.
Then $X_2^2\leq2Q+8u^2$ and
$\sum_{i\geq3}X_i^2\leq2P+4Q+16u^2$, since
$\sum_{i\geq3}(x_i/x_2)^2\leq\epsilon^2\leq1$.
Thus $|\mathbf{X}_{\perp}|^2+u^2\leq25(P+Q+u^2)$, whereas
\eqref{eq:squares} is at least
$2\epsilon^2x_2(P+Q+u^2)$.
The constants in these estimates are uniform as the regular boundary
approaches $I$, including its endpoints.

Write $g_{ij}$ for the coefficients of $\Phi^*g$, and let $g_\rho$
have coefficients $g_{ij}(\rho \mathbf{x})$. The pullback under
$\mathbf{x}\mapsto\Phi(\rho \mathbf{x})$ is $\rho^2g_\rho$; the constant factor does
not change the connection or Hessian. After decreasing \(\rho\) sufficiently, assume \(0<\rho\le1\) and \(\rho|\mathbf y|_{\mathrm E}\le1\) throughout \(C_\epsilon\). All constants below are
uniform in the boundary point and in small $\rho$; the letter $C$
may denote different positive constants.

\emph{The metric along the axis and its transverse expansion.}
%At $(t,\mathbf{0})$ the coordinate vectors are mapped by $d\Phi$ to
%$E_1(t),\ldots,E_n(t)$: this follows from $\Phi(t,\mathbf{0})=\gamma(t)$
%and $d(\exp_{\gamma(t)})_0=\operatorname{Id}$. Since the chosen
%parallel frame is orthonormal,
%\[
%g_{ij}(t,\mathbf{0})
%=g_{\gamma(t)}\bigl(E_i(t),E_j(t)\bigr)=\delta_{ij}.
%\]
%Thus this equality holds along the axis by construction; it is
%not an assumption on the metric away from the axis.
%Parallelness gives $\Gamma^k_{1i}(t,\mathbf{0})=0$, and symmetry gives
%$\Gamma^k_{i1}(t,\mathbf{0})=0$. For every fixed normal vector $v$, the
%curve $s\mapsto\Phi(t,sv)$ is a geodesic, so its equation at $s=0$
%gives $\sum_{\alpha,\beta\geq2}
%\Gamma^k_{\alpha\beta}(t,\mathbf{0})v_\alpha v_\beta=0$.
%Since this holds for every $v$, all the remaining Christoffel
%symbols vanish on the axis as well.
When restricted to $\mathbf{y}=\mathbf{0}$, $d\Phi(\partial_i)=E_i$, so $g_{ij}(t,\mathbf{0})=\delta_{ij}$. Since the frame is parallel, $\Gamma^k_{1i}(t,\mathbf{0})=0$. Moreover, each normal ray $s\mapsto\Phi(t,sv)$ is a geodesic; hence
\[
\sum_{i,j=2}^n\Gamma^k_{ij}(t,\mathbf{0})v_iv_j=0
\]
for every normal vector $v$. By symmetry in $i,j$, all Christoffel symbols vanish on the axis.

The normal exponential construction gives
\[
g_{11}(t,\mathbf{y})=1-\langle A(t)\mathbf{y},\mathbf{y}\rangle+O(|\mathbf{y}|^3),
\qquad g_{1{i}}(t,\mathbf{y})=O(|\mathbf{y}|^2)\quad({i}\geq2).
\]
The first identity follows by varying the base point
of a normal geodesic: the resulting Jacobi field has initial data
$E_1(t),0$. Taylor expansion in the transverse variable, followed
by the homothety, therefore yields
\begin{equation}
\label{eq:christoffel}
\begin{gathered}
\Gamma^{i}_{11}(g_\rho)
=\rho^2(A(\rho t)\mathbf{y})^{i}+O(\rho^3|\mathbf{y}|^2),\qquad {i}\geq2,\\
\Gamma^1_{11}(g_\rho)=O(\rho^3|\mathbf{y}|^2),\qquad
|\Gamma^k_{ij}(g_\rho)|\leq C\rho^2|\mathbf{y}|.
\end{gathered}
\end{equation}
We estimate the Riemannian Hessian $
\operatorname{Hess}_{g_\rho}F
=
\bigl(\mathbf x\mapsto\Phi(\rho\mathbf x)\bigr)^*
\operatorname{Hess}_g
\left(F\circ(\rho^{-1}\Phi^{-1})\right) $. \footnote{Explicitly, $ \operatorname{Hess}_{g_\rho}F\big|_{\mathbf x}
	(\mathbf X,\mathbf Y)
	=
	\operatorname{Hess}_g
	\left(F\circ(\rho^{-1}\Phi^{-1})\right)
	\big|_{\Phi(\rho\mathbf x)}\left(
	\rho\,d\Phi_{\rho\mathbf x}\mathbf X,\,
	\rho\,d\Phi_{\rho\mathbf x}\mathbf Y
	\right) $.} 

For a tangent vector $\mathbf{X}=X_1\partial_1+\mathbf{X}_{\perp}$, we have
\begin{equation}
	\label{eq:hessian-initial}
	\begin{aligned}
		\operatorname{Hess}_{g_\rho}F(\mathbf{X},\mathbf{X})
		& = \langle D_{\mathbf X}\nabla F,\mathbf X\rangle_{\mathrm E}-(D_{\mathbf X}g_\rho)
		\bigl(\mathbf X,\nabla^{g_\rho}F\bigr)
		+\frac12(D_{\nabla^{g_\rho}F}g_\rho)
		\bigl(\mathbf X,\mathbf X\bigr)
	\end{aligned}
\end{equation}
Here \(D\) is the Euclidean connection, \(\nabla F\) is the Euclidean gradient, and \(\nabla^{g_\rho}F\) is the gradient for \(g_\rho\). The smooth Fermi expansions give
\[
\begin{aligned}
	\|g_\rho-\langle\cdot,\cdot\rangle_{\mathrm E}\|
	&\le C\rho^2|\mathbf y|_{\mathrm E}^2,\\
	\|D_{e_1}g_\rho\|
	&\le C\rho^3|\mathbf y|_{\mathrm E}^2,\\
	\|D_{\mathbf X_\perp}g_\rho\|
	&\le C\rho^2|\mathbf y|_{\mathrm E}|\mathbf X_\perp|_{\mathrm E}.
\end{aligned}
\]
Moreover,
\[
|\nabla^{g_\rho}F|_{\mathrm E}\le C|\nabla F|_{\mathrm E},
\qquad
|\nabla^{g_\rho}F-\nabla F|_{\mathrm E}
\le C\rho^2|\mathbf y|_{\mathrm E}^2|\nabla F|_{\mathrm E}.
\]
Since $\mathbf X=\langle\mathbf X,e_1\rangle_{\mathrm E}e_1+\mathbf X_\perp$,
the second term in \eqref{eq:hessian-initial} satisfies
\[
\begin{aligned}
\left|(D_{\mathbf X}g_\rho)
(\mathbf X,\nabla^{g_\rho}F)\right|
&\le C|\nabla F|_{\mathrm E}|\mathbf X|_{\mathrm E}
\left(
\rho^3|\mathbf y|_{\mathrm E}^2
|\langle\mathbf X,e_1\rangle_{\mathrm E}|
+\rho^2|\mathbf y|_{\mathrm E}|\mathbf X_\perp|_{\mathrm E}
\right)\\
&\le C|\nabla F|_{\mathrm E}
\Bigl[
\rho^3|\mathbf y|_{\mathrm E}^2
\langle\mathbf X,e_1\rangle_{\mathrm E}^2+
\rho^2|\mathbf y|_{\mathrm E}
\left(
|\langle\mathbf X,e_1\rangle_{\mathrm E}|
|\mathbf X_\perp|_{\mathrm E}
+|\mathbf X_\perp|_{\mathrm E}^2
\right)
\Bigr],
\end{aligned}
\]
where we used $\rho|\mathbf y|_{\mathrm E}\le1$.

For the third term in \eqref{eq:hessian-initial}, bilinearity gives
\[
\begin{aligned}
	\frac12(D_{\nabla^{g_\rho}F}g_\rho)(\mathbf X,\mathbf X)
	&=\frac12\langle\mathbf X,e_1\rangle_{\mathrm E}^2
	(D_{\nabla^{g_\rho}F}g_\rho)(e_1,e_1)\\
	&\quad+\langle\mathbf X,e_1\rangle_{\mathrm E}
	(D_{\nabla^{g_\rho}F}g_\rho)(e_1,\mathbf X_\perp)\\
	&\quad+\frac12(D_{\nabla^{g_\rho}F}g_\rho)
	(\mathbf X_\perp,\mathbf X_\perp),
\end{aligned}
\]
where the last two terms are bounded in absolute value by $ C\rho^2|\mathbf y|_{\mathrm E}|\nabla F|_{\mathrm E}
\left(
|\langle\mathbf X,e_1\rangle_{\mathrm E}|
|\mathbf X_\perp|_{\mathrm E}
+|\mathbf X_\perp|_{\mathrm E}^2
\right)$. For the first term $\frac12\langle\mathbf X,e_1\rangle_{\mathrm E}^2
(D_{\nabla^{g_\rho}F}g_\rho)(e_1,e_1)$, differentiate
\[
g_\rho(e_1,e_1)
=1-\rho^2\langle A(\rho t)\mathbf{y},\mathbf{y}\rangle_{\mathrm E}
+O(\rho^3|\mathbf y|_{\mathrm E}^3).
\]
The normal derivative of the quadratic expression gives
\[
\frac12(D_{\nabla^{g_\rho}F}g_\rho)(e_1,e_1)
=-\rho^2
\langle\nabla^{g_\rho}F,(0,A(\rho t)\mathbf{y})\rangle_{\mathrm E}
+O\!\left(
\rho^3|\mathbf y|_{\mathrm E}^2|\nabla F|_{\mathrm E}
\right).
\]
Differentiating $A(\rho t)$ contributes an additional factor $\rho$,
so that contribution is included in the remainder. Replacing
$\nabla^{g_\rho}F$ by $\nabla F$ costs at most
\[
\begin{aligned}
	\rho^2\left|
	\langle\nabla^{g_\rho}F-\nabla F,(0,A(\rho t)\mathbf{y})\rangle_{\mathrm E}
	\right|
	&\le C\rho^4|\mathbf y|_{\mathrm E}^3|\nabla F|_{\mathrm E}\le C\rho^3|\mathbf y|_{\mathrm E}^2|\nabla F|_{\mathrm E}.
\end{aligned}
\]
Thus the sum of the second and third terms in \eqref{eq:hessian-initial} equals
\[
-\rho^2\langle\mathbf X,e_1\rangle_{\mathrm E}^2
\langle\nabla F,(0,A(\rho t)\mathbf{y})\rangle_{\mathrm E}
\]
up to a remainder having absolute value at most
\[
\begin{aligned}
	C|\nabla F|_{\mathrm E}\Bigl[
	&\rho^3|\mathbf y|_{\mathrm E}^2
	\langle\mathbf X,e_1\rangle_{\mathrm E}^2+\rho^2|\mathbf y|_{\mathrm E}
	\left(
	|\langle\mathbf X,e_1\rangle_{\mathrm E}|
	|\mathbf X_\perp|_{\mathrm E}
	+|\mathbf X_\perp|_{\mathrm E}^2
	\right)
	\Bigr].
\end{aligned}
\]

On the regular boundary,
$|\nabla F|_{\mathrm E}\le4\epsilon x_2$ and
$|\mathbf y|_{\mathrm E}\le\sqrt2x_2$.
Consequently, this remainder is bounded by
\[
\begin{aligned}
	C\epsilon\rho^2x_2^2\Bigl[
	&\rho x_2\langle\mathbf X,e_1\rangle_{\mathrm E}^2+
	|\langle\mathbf X,e_1\rangle_{\mathrm E}|
	|\mathbf X_\perp|_{\mathrm E}
	+|\mathbf X_\perp|_{\mathrm E}^2
	\Bigr].
\end{aligned}
\]
Using $0<\rho\le1$ and $x_2\le\theta\le1$, we have
\[
\begin{aligned}
	\rho x_2^3\langle\mathbf X,e_1\rangle_{\mathrm E}^2
	&\le x_2\theta^2\langle\mathbf X,e_1\rangle_{\mathrm E}^2,\\
	x_2^2|\langle\mathbf X,e_1\rangle_{\mathrm E}|
	|\mathbf X_\perp|_{\mathrm E}
	&\le\frac{x_2}{2}
	\left(
	\theta^2\langle\mathbf X,e_1\rangle_{\mathrm E}^2
	+|\mathbf X_\perp|_{\mathrm E}^2
	\right),\\
	x_2^2|\mathbf X_\perp|_{\mathrm E}^2
	&\le x_2|\mathbf X_\perp|_{\mathrm E}^2.
\end{aligned}
\]

Together with~\eqref{eq:model-identities}, these bounds give
\begin{equation}
\label{eq:hessian-perturbation}
\begin{aligned}
\operatorname{Hess}_{g_\rho}F(\mathbf{X},\mathbf{X})
& \geq D^2F(\mathbf{X},\mathbf{X})-\rho^2\langle\mathbf X,e_1\rangle_{\mathrm E}^2
\langle\nabla F,(0,A(\rho t)\mathbf{y})\rangle_{\mathrm E}-C\epsilon\rho^2x_2
\left(
|\mathbf X_\perp|_{\mathrm E}^2
+\theta^2\langle\mathbf X,e_1\rangle_{\mathrm E}^2\right).
\end{aligned}
\end{equation}
%Indeed, the terms involving $\mathbf{X}_{\perp}$ are bounded by
%$C\epsilon\rho^2x_2^2(|X_1||\mathbf{X}_{\perp}|+|\mathbf{X}_{\perp}|^2)$; the remaining error involving $ \langle\mathbf X,e_1\rangle_{\mathrm E} $
%is bounded by $C\epsilon\rho^3x_2^3X_1^2$. Use $x_2\leq\theta\leq1$
%and $2\theta|X_1||\mathbf{X}_{\perp}|\leq\theta^2X_1^2+|\mathbf{X}_{\perp}|^2$ to obtain the stated
%bound. 

%It controls directions close to the axis and transverse directions simultaneously.

\emph{Case 1: constant sectional curvature.}
If $A(t)=\lambda\operatorname{Id}$, then
\[
-\langle\nabla F,(0,A(\rho t)\mathbf{y})\rangle_{\mathrm E} = -\lambda\epsilon^2x_2^2\theta
\geq-|\lambda|\epsilon^2x_2\theta^2.
\]
Equations~\eqref{eq:lower-bound} and~\eqref{eq:hessian-perturbation}
imply
\begin{equation}
\label{eq:constant-hessian}
\operatorname{Hess}_{g_\rho}F(\mathbf{X},\mathbf{X})
\geq x_2\left(\frac{2\epsilon^2}{25}-C\epsilon\rho^2\right)
\bigl(|\mathbf{X}_{\perp}|^2+\theta^2X_1^2\bigr)>0
\end{equation}
for every nonzero boundary tangent vector $\mathbf{X}$ when $\rho$ is
sufficiently small.

\emph{Case 2: a simple maximal eigenvalue.}
Let $\delta=\lambda_2-\max_{i\geq3}\lambda_i>0$.
The identities in~\eqref{eq:model-identities} give
\[
\begin{aligned}
-\langle\nabla F,(0,A(0)\mathbf{y})\rangle_{\mathrm E}
&=2\sum_{i=3}^n(\lambda_2-\lambda_i)x_i^2
-\lambda_2\epsilon^2x_2^2\theta\\
&\geq2\delta\epsilon^2x_2^2(1-\theta)
-|\lambda_2|\epsilon^2x_2^2\theta\\
&\geq2\delta\epsilon^2x_2^2
-(2\delta+|\lambda_2|)\epsilon^2x_2\theta^2.
\end{aligned}
\]
Moreover, $\|A(\rho t)-A(0)\|_{\mathrm{op}}\leq C\rho$, so
\[
|\langle\nabla F,(0,A(\rho t)\mathbf{y}- A(0)\mathbf{y})\rangle_{\mathrm E}|\leq C\epsilon\rho x_2^2.
\]
Choose $\rho$ so that $C\rho\leq\delta\epsilon$. Substitution in
\eqref{eq:hessian-perturbation} yields
\begin{equation}
\label{eq:gap-hessian}
\begin{aligned}
\operatorname{Hess}_{g_\rho}F(\mathbf{X},\mathbf{X})
&\geq x_2\left(\frac{2\epsilon^2}{25}-C\epsilon\rho^2\right)
\bigl(|\mathbf{X}_{\perp}|^2+\theta^2X_1^2\bigr)+\delta\epsilon^2\rho^2x_2^2X_1^2>0
\end{aligned}
\end{equation}
for every nonzero boundary tangent vector $\mathbf{X}$, after a further
decrease of $\rho$. The same positive gap \(\delta\) used in Proposition~\ref{prop:jacobi} also supplies the positive axial term in the Hessian estimate.

In both cases, the same sufficiently small $\rho$ works on the
entire regular boundary. Near each regular boundary point, the body is given by \(F\circ(\rho^{-1}\Phi^{-1})\le0\), with nonvanishing differential. Therefore positivity of its Hessian on boundary tangent vectors implies positive definiteness of the second fundamental form. Strict Euclidean convexity alone gives stability on compact subsets
away from $I$; the estimates above also control approach to $I$
and its endpoints.

\emph{Local convexity at the singular segment.} Since
$C_\epsilon\subset\RR e_1\oplus \mathop{\mathrm{Cone}}\nolimits_\epsilon$, the set
$T_q\mathcal C_{\epsilon,\rho}\cap(-T_q\mathcal C_{\epsilon,\rho})$ at a point of the segment
is contained in $\RR\dot\gamma$. At either endpoint it is zero:
the additional supporting inequality $-1\leq x_1\leq1$ excludes
one of the two directions tangent to the geodesic segment. 

%These statements follow by applying
%the derivative of $x\mapsto\Phi(\rho x)$ to the model tangent cones.

Fix a singular point $q$ and choose a small strongly convex ball $B$ about $q$ whose preimage under $\Phi$ is Euclidean convex. Then $B\cap\operatorname{int}\mathcal C_{\epsilon,\rho}$ is path connected and dense in $B\cap\mathcal C_{\epsilon,\rho}$. Fix one endpoint in this interior and move the other along an interior path starting there. By continuity, a first failure of the joining geodesics to remain in the interior would give a geodesic contained in $\mathcal C_{\epsilon,\rho}$ and touching its boundary between its endpoints. At a regular contact, $F\circ(\rho^{-1}\Phi^{-1})$ has a local maximum along the geodesic, contradicting its positive tangential Hessian. At a singular contact, both velocity directions belong to the tangent cone. This excludes contact at the segment endpoints and forces the velocity to be parallel to $\dot\gamma$ elsewhere. Uniqueness then forces the geodesic to coincide with $\gamma$, contradicting that its endpoints lie in \(\operatorname{int}\mathcal C_{\epsilon,\rho}\). Hence all joining geodesics remain in the interior, and approximation gives convexity of $B\cap\mathcal C_{\epsilon,\rho}$. At regular boundary points, local convexity follows from the positive-definite second fundamental form.

%Fix a singular point $q$. Choose a small strongly convex ball $B$ about $q$ whose coordinate preimage is Euclidean convex; this is possible because the coordinate Hessian of $d(q,\cdot)^2$ is positive definite near $q$. Thus $B\cap\operatorname{int}\mathcal C$ is path connected and dense in $B\cap\mathcal C$. Fix $a$ in this interior and move the other endpoint along an interior path from $a$. The joining minimizing geodesics vary continuously. If one first fails to remain in the interior, it lies in $\mathcal C$ and touches its boundary between its endpoints.
%
%At a regular contact, \(F\circ(\rho^{-1}\Phi^{-1})\) has a local maximum along the geodesic, contradicting its positive tangential Hessian. At a singular contact, both velocity directions lie in the tangent cone. At an endpoint this would force the velocity to vanish, which is impossible for a nonconstant geodesic; at any other point of the segment, it forces the velocity to be parallel to \(\dot\gamma\). By uniqueness, the geodesic then coincides with $\gamma$, contradicting its interior endpoints. Hence all joining geodesics remain in the interior. Approximation gives convexity of $B\cap\mathcal C$. At regular boundary points, local convexity follows from the positive-definite second fundamental form.
\end{proof}

\section{Proof of Theorem~\ref{thm:main}}
\begin{proof}
%Choose $0<\epsilon\leq1$ as small as desired. In the first case
%any parallel normal frame is suitable; in the second case
%choose the frame of Proposition~\ref{prop:jacobi}. The corresponding parallel
%cones are Jacobi convex on a sufficiently short interval by
%Section~2. Proposition~\ref{prop:local} supplies the locally convex
%body $\mathcal C=\Phi(\rho C_\epsilon)$ for all sufficiently small
%$\rho$. 
Fix $0<\epsilon\leq1$ and choose the parallel frame as above. Proposition~\ref{prop:local} gives a locally convex body $\mathcal C_{\epsilon,\rho}=\Phi(\rho C_\epsilon)$ for all sufficiently small $\rho$. We choose $\rho$ smaller to obtain global convexity.

%Choose a relatively compact neighborhood $W$ of $p$, inside the
%prescribed neighborhood, and $i_0>0$ such that
%$\operatorname{inj}_M(q)>i_0$ for every $q\in W$. This uniform
%normal-neighborhood property is local and does not require
%completeness. On a compact coordinate neighborhood $Q\subset U$ of $0$
%with $\Phi(Q)\subset W$, choose $L>0$ such that
%$|d\Phi_z(v)|_g\leq L|v|$ for $z\in Q$. Choose $\rho_0>0$ so that
%$\rho_0C_\epsilon\subset Q$ and Proposition~\ref{prop:local} applies
%at every scale $0<\rho\leq\rho_0$. Set
%\[
%s=\min\left\{1,
%\frac{i_0}{2L\rho_0\operatorname{diam}_{\mathrm E}(C_\epsilon)}\right\}
%\in(0,1],\qquad \rho=s\rho_0>0.
%\]
%Thus the additional homothety has the strictly positive ratio $s$,
%and the final body is
%\[
%\mathcal C=\Phi(s\rho_0 C_\epsilon).
%\]
%Since $C_\epsilon$ is convex and contains $0$, we still have
%$s\rho_0C_\epsilon\subset\rho_0C_\epsilon\subset Q$,
%and local convexity holds at this scale. Its scale satisfies
%\[
%L(s\rho_0)\operatorname{diam}_{\mathrm E}(C_\epsilon)
%\leq\frac{i_0}{2}<i_0.
%\]
%For $a,b\in C_\epsilon$, the image of the straight segment
%$\rho[(1-t)a+tb]$, $0\leq t\leq1$, lies in $\mathcal C$ and has length at most
%$L\rho|a-b|<i_0$.
Choose a relatively compact neighborhood $W$ of $p$ inside the prescribed neighborhood, and $i_0>0$ such that $\operatorname{inj}_M(q)>i_0$ on $W$. Fix a compact coordinate neighborhood $Q$ of $0$ with $\Phi(Q)\subset W$, and let $L$ bound $|d\Phi|$ on $Q$. Choose $\rho>0$ so small that Proposition~\ref{prop:local} applies, $\rho C_\epsilon\subset Q$, and
\[
L\rho\,\operatorname{diam}_{\mathrm E}C_\epsilon<i_0.
\]
Images of straight segments then join any two points of $\mathcal C_{\epsilon,\rho}$ inside $\mathcal C_{\epsilon,\rho}$ with length less than $i_0$.

%Compactness of $\mathcal C$, the Arzel\`a--Ascoli theorem, and
%lower semicontinuity of length give an intrinsic shortest curve
%between any two of its points, of length less than $i_0$.
%Every sufficiently short subarc can be replaced, in a locally
%convex neighborhood, by the ambient minimizing geodesic with the
%same endpoints. Minimality of the curve forces each such subarc
%to be that geodesic. Hence the intrinsic shortest curve, parametrized
%by arc length, is an ambient geodesic. Its length is smaller than
%the injectivity radius at its initial point, so it minimizes length
%in $M$. Any ambient minimizing geodesic with these endpoints has
%the same length, less than $i_0$. Injectivity of the exponential map
%on this radius forces it to coincide with the intrinsic shortest
%curve. Thus every ambient minimizing geodesic between the two points
%is contained in $\mathcal C$, proving the required convexity.
Compactness and Arzelà–Ascoli give an intrinsic shortest curve between any two points of $\mathcal C_{\epsilon,\rho}$. Its length is less than $i_0$. Local convexity makes each sufficiently short subarc an ambient minimizing geodesic. Thus the curve is an ambient geodesic; since its length is smaller than the injectivity radius at its initial point, it is the unique ambient minimizing geodesic.

Finally, $\Phi(\rho I)=\gamma([-\rho,\rho])$ lies in $\partial\mathcal C_{\epsilon,\rho}$.
At every point in the relative interior of this
segment, the differential of $\Phi$ at $(t,0,...,0)$ identifies the coordinate frame with the parallel orthonormal frame. Therefore
\[
T_{\gamma(t)}\mathcal C_{\epsilon,\rho}
=\RR\dot\gamma(t)\oplus K_t,\qquad |t|<\rho.
\]
Here \(K_t\) is the parallel translate to \(\dot\gamma(t)^\perp\) of the cone \(\mathop{\mathrm{Cone}}\nolimits_\epsilon\) at \(\gamma(0)\). 
%It remains to apply Propositions~\ref{prop:closed-convex-hypersurface-skinny-cones}and~\ref{prop:local}. 
\end{proof}

\section{Proof of Theorem~\ref{thm:four}}
\begin{proof}
	Applying Theorem~\ref{thm:main}, we may assume that no normal Jacobi operator has a simple maximal eigenvalue. If all normal Jacobi operators are
	scalar, the argument of Corollary~\ref{cor:three} gives constant
	sectional curvature on $M$ and Theorem~\ref{thm:main} applies. Otherwise, choose a short geodesic along
	which the maximal eigenvalue is double and separated from the simple
	minimal eigenvalue. In the parallel identification of the normal
	spaces with $\RR^3$, write
	\[
	A(t)V=\lambda_{\max}(t)V-(\lambda_{\max}(t)-\lambda_{\min}(t))
	\langle V,w(t)\rangle w(t),\qquad \lambda_{\max}(t)>\lambda_{\min}(t).
	\]
	Here $\lambda_{\max}$ and $\lambda_{\min}$ are smooth eigenvalues of $A$ and $w$ is a smooth unit $\lambda_{\min}$-eigenvector.
	
	Denote by $\pi_t$ the orthogonal projection onto $w(t)^\perp$. We construct, on some interval $|t|\le b/2$, a skinny cone $K$
	preserved by $\pi_t$.
	If $w$ is parallel on a segment, a circular
	cone about a parallel maximal eigenvector suffices. Otherwise recenter at a point
	where $\dot w(0)\ne0$, choose $E_4(0)=w(0)$ and
	$E_3(0)=\dot w(0)/|\dot w(0)|$, and, for small $\alpha>0$, set
	\[
	v_\pm(0)=\cos\alpha\,E_2(0)\pm\sin\alpha\,E_3(0),
	\qquad
	\nu_\pm(0)=\sin\alpha\,E_2(0)\mp\cos\alpha\,E_3(0).
	\]
	Extend these vectors by parallel transport along the spherical curve $w$:
	\[
	\dot v_\pm=-\langle v_\pm,\dot w\rangle w,
	\qquad
	\dot\nu_\pm=-\langle\nu_\pm,\dot w\rangle w.
	\]
	Parallel transport preserves inner products, so each triple $v_\pm,\nu_\pm,w$ remains orthonormal. After shortening the interval,
	$\langle v_\pm(s),\dot w(t)\rangle$ and
	$\langle\nu_\pm(s),\dot w(t)\rangle$ retain opposite signs,
	uniformly in $s,t$.
	
	Let $K$ be the closed convex cone generated by $v_\pm(s)$, $|s|\le b$,
	for sufficiently small $b>0$ (Figure~\ref{fig:cone K}). We show that, for $|s|<b$, the vectors $\nu_\pm(s)$ are the inward unit normals to $\partial K$ along the rays $\mathbb R_{>0}v_\pm(s)$, respectively. Indeed, differentiating twice gives
	\[
	\langle\nu_\pm(s),v_\pm(r)\rangle>0\quad(r\ne s),
	\qquad
	\langle\nu_\pm(s),v_\mp(r)\rangle>0,
	\]
	since the same-arc second derivative at $r=s$ is
	$-\langle v_\pm,\dot w\rangle\langle\nu_\pm,\dot w\rangle>0$.
	For sufficiently small \(b\), continuity makes the same-arc pairings strictly convex in \(r\) and keeps the cross-arc pairings positive.
	\begin{figure}[htbp]
		\centering
		\includegraphics[width=0.5\textwidth]{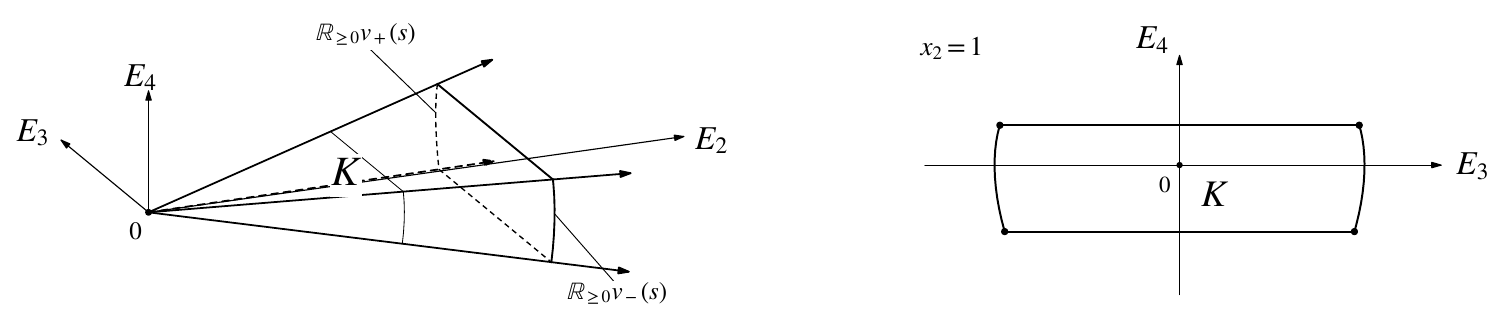}
		\caption{The cone $K$}
		\label{fig:cone K}
	\end{figure}

	Let $\nu_{\mathrm{top}}$ be the inward unit normal to the plane spanned by $v_+(-b)$ and $v_-(b)$, chosen so that $\nu_{\mathrm{top}}\to-w(0)$ as \(b\downarrow0\). The equations defining $v_\pm$ and $\nu_\pm$ imply that $\langle\nu_{\mathrm{top}},v_+(s)\rangle$ is strictly increasing and $\langle\nu_{\mathrm{top}},v_-(s)\rangle$ is strictly decreasing. Since these functions vanish at $-b$ and $b$, respectively, the plane supports both arcs. The other crossed chord is treated similarly. Away from the vertex, $\partial K$ consists of the two conical surfaces generated by $v_\pm(s)$, $|s|\le b$, and the two planar faces generated by the crossed endpoint pairs $v_+(-b),v_-(b)$ and $v_+(b),v_-(-b)$. Thus $K$ is full-dimensional,
	pointed, and arbitrarily skinny.
	
	For $|t|\le b/2$,
	\[
	\langle A(t)v_\pm(s),\nu_\pm(s)\rangle
	=-(\lambda_{\max}(t)-\lambda_{\min}(t))
	\langle v_\pm(s),w(t)\rangle
	\langle\nu_\pm(s),w(t)\rangle\ge0.
	\]
	The last two factors vanish at $t=s$ and are integrals of functions
	with opposite signs. On each chord face the same inward inequality
	is strict at nonzero points: its inward normal has scalar product with $w(t)$ of the
	sign opposite to both endpoint scalar products. At each junction ray, every supporting normal is a nonnegative combination of the adjacent arc and chord normals. Hence every supporting face
	satisfies~\eqref{eq:inward}. The projection argument of
	Lemma~\ref{lem:jacobi}, applied to the linear flow, shows that
	$\exp(r(A(t)-\lambda_{\max}(t)\operatorname{Id}))K\subset K$ for $r\ge0$.
	Taking $r\to\infty$ gives $\pi_tK\subset K$. Consequently, for a
	uniform $c\ge0$,
	\begin{equation}
		\label{eq:four-positive}
		(A(t)+c\operatorname{Id})K\subset K,
	\end{equation}
	because $A(t)=\lambda_{\min}(t)\operatorname{Id}+(\lambda_{\max}(t)-\lambda_{\min}(t))\pi_t$.
	
	We now construct a convex body whose transverse tangent cones along $\gamma$ are parallel copies of $K$. Consider the Fermi coordinates $(t,\mathbf{y})$ as in \eqref{eq:fermi} with $n=4$. Let $a\in (0,\min\{1,b/2\})$ be chosen sufficiently small below. For $\eta>0$ and every unit $\nu\in K^*$, define
	\[
	F_{\nu,\eta}(t,\mathbf{y})=\frac{|\mathbf{y}|^2+\eta}{a^2-t^2}-\langle\nu,\mathbf{y}\rangle,
	\qquad
	C_\eta=\bigcap_{\nu\in K^*,\ |\nu|=1}
	\{(t,\mathbf{y}):|t|<a,\ F_{\nu,\eta}(t,\mathbf{y})\le0\}.
	\]
	
	At a point $(t,\mathbf{y})\in\partial C_\eta$, choose a unit $\nu\in K^*$ such that $F_{\nu,\eta}(t,\mathbf{y})=0$. Then
	\[
	\frac{|\mathbf{y}|^2+\eta}{a^2-t^2}=\langle\nu,\mathbf{y}\rangle\le |\mathbf{y}|,
	\qquad |\mathbf{y}|\le a^2-t^2,\qquad
	|\partial_tF_{\nu,\eta}|\le2a,\quad |d_\mathbf{y}F_{\nu,\eta}|\le3.
	\]
	Applying~\eqref{eq:four-positive},
	$\langle A(t)\mathbf{y},\nu\rangle\ge-c(|\mathbf{y}|^2+\eta)/(a^2-t^2)$.
	Applying~\eqref{eq:hessian-initial} with \(g_\rho\) replaced by \(\Phi^*g\) and \(F\) by \(F_{\nu,\eta}\), we have 
	\[
	\begin{aligned}
		\operatorname{Hess}_{\Phi^*g}F_{\nu,\eta}(\mathbf X,\mathbf X)
		&=D^2F_{\nu,\eta}(\mathbf X,\mathbf X)-(D_{\mathbf X}(\Phi^*g))
		\bigl(\mathbf X,\nabla^{\Phi^*g}F_{\nu,\eta}\bigr)+\frac12(D_{\nabla^{\Phi^*g}F_{\nu,\eta}}(\Phi^*g))
		\bigl(\mathbf X,\mathbf X\bigr).
	\end{aligned}
	\]
	Here, for $\mathbf{X}=X_1\partial_1+\mathbf{X}_{\perp}$, we have 
	\[
	\begin{aligned}
		D^2F_{\nu,\eta}(\mathbf X,\mathbf X)
		&=\frac{2}{a^2-t^2}
		\left|\mathbf X_\perp+\frac{2tX_1}{a^2-t^2}\mathbf y\right|^2\quad+\frac{2(|\mathbf y|^2+\eta)}{(a^2-t^2)^2}X_1^2
		+\frac{8t^2\eta}{(a^2-t^2)^3}X_1^2,
	\end{aligned}
	\]
	\[
	\begin{aligned}
		\left|
		(D_{\mathbf X}(\Phi^*g))
		\bigl(\mathbf X,\nabla^{\Phi^*g}F_{\nu,\eta}\bigr)
		\right|&\le
		\left\|D_{\mathbf X}(\Phi^*g)\right\|_{\mathrm{op}}
		|\mathbf X|_{\mathrm E}
		\left|\nabla^{\Phi^*g}F_{\nu,\eta}\right|_{\mathrm E}\\
		&\le
		C\left(|\mathbf y|^2|X_1|
		+|\mathbf y||\mathbf X_\perp|\right)
		\left(|X_1|+|\mathbf X_\perp|\right)\\
		&\le C\left(
		|\mathbf y|^2X_1^2
		+|\mathbf y||X_1||\mathbf X_\perp|
		+|\mathbf y||\mathbf X_\perp|^2
		\right),
	\end{aligned}
	\]
	and
	\[
	\begin{aligned}
		&\frac12(D_{\nabla^{\Phi^*g}F_{\nu,\eta}}(\Phi^*g))
		(\mathbf X,\mathbf X)=
		-X_1^2\left\langle\nabla F_{\nu,\eta},
		(0,A(t)\mathbf y)\right\rangle_{\mathrm E}+
		O\!\left(
		|\mathbf y|^2X_1^2
		+|\mathbf y||X_1||\mathbf X_\perp|
		+|\mathbf y||\mathbf X_\perp|^2
		\right).
	\end{aligned}
	\]
	Hence,
	\[
	\begin{aligned}
		\operatorname{Hess}_{\Phi^*g}F_{\nu,\eta}(\mathbf{X},\mathbf{X})
		&\ge D^2F_{\nu,\eta}(\mathbf X,\mathbf X)
		-X_1^2\left\langle
		\nabla F_{\nu,\eta},(0,A(t)\mathbf y)
		\right\rangle_{\mathrm E}-C\left(
		|\mathbf y|^2X_1^2
		+|\mathbf y||X_1||\mathbf X_\perp|
		+|\mathbf y||\mathbf X_\perp|^2
		\right)\\
		&\ge \frac{2}{a^2-t^2}\left|\mathbf{X}_{\perp}+\frac{2tX_1}{a^2-t^2}\mathbf{y}\right|^2
		+\frac{2(|\mathbf{y}|^2+\eta)}{(a^2-t^2)^2}X_1^2-C\left(\frac{|\mathbf{y}|^2+\eta}{a^2-t^2}X_1^2
		+|\mathbf{y}||X_1||\mathbf{X}_{\perp}|+|\mathbf{y}||\mathbf{X}_{\perp}|^2\right).
	\end{aligned}
	\]
	Since $|\mathbf{y}|\le a^2-t^2\leq a^2$ and $|t|\le a$, completing the square absorbs the error terms for sufficiently small $a$ and yields
	\[
	\operatorname{Hess}_{\Phi^*g}F_{\nu,\eta}(\mathbf{X},\mathbf{X})
	\ge \frac{1}{a^2-t^2}\left|\mathbf{X}_{\perp}+\frac{2tX_1}{a^2-t^2}\mathbf{y}\right|^2
	+\frac{|\mathbf{y}|^2+\eta}{(a^2-t^2)^2}X_1^2>0
	\qquad(\mathbf{X}\ne0).
	\]
	At each boundary point, compactness gives a common neighborhood where the defining functions that vanish there, together with those for nearby normals, have positive-definite Hessians. The remaining inequalities stay strict. Hence $\Phi(C_\eta)$ is locally convex. Every $C_\eta$ is Euclidean convex by the same
	square completion, and compact since $\eta\le (a^2-t^2)^2/4$ keeps it away
	from $|t|=a$. It has nonempty interior for sufficiently small $\eta$.
	For sufficiently small \(a\), compactness and Euclidean convexity of \(C_\eta\) give a shortest curve inside \(\Phi(C_\eta)\) between any two points, with length below the injectivity radius. Local convexity makes this curve an ambient geodesic, hence the unique ambient minimizing geodesic. Thus \(\Phi(C_\eta)\) is convex.
	
	Define $C_0$ by \[
	C_0=\left\{(t,\mathbf y):|t|\le a,\quad
	|\mathbf{y}|^2\le(a^2-t^2)\langle\nu,\mathbf y\rangle
	\text{ for every unit }\nu\in K^*\right\}.
	\] 
	Radial contraction in $\mathbf{y}$ shows that
	$C_0=\overline{\bigcup_{\eta>0}C_\eta}$, with the axis obtained by
	approximation from the interior of $K$.
	Hence $\Phi(C_0)$ is a compact convex body with nonempty
	interior. With \(\eta=0\), the same Hessian estimate remains positive at boundary points with \(\mathbf{y}\ne\mathbf{0}\). Along a boundary geodesic through such a point, some \(F_{\nu,0}\) would attain a local maximum \(0\), contradicting the positive Hessian estimate. Thus every nontrivial boundary geodesic segment lies on the axis. Since $C_0\subset\mathbb R\times K$, its axis segment lies in its boundary, and $T_{(t,\mathbf{0})}C_0\subset\mathbb R\times K$ for $|t|<a$. Conversely, for $v\in\operatorname{int}K$ and any $X_1$, the points $(t+sX_1,sv)$ belong to $C_0$ for sufficiently small $s>0$. Taking closures gives the reverse inclusion. Applying $d\Phi_{(t,\mathbf{0})}$, we obtain
	\[
	T_{\gamma(t)}\Phi(C_0)
	=\RR\dot\gamma(t)\oplus \mathcal{P}_{0\to t}K,
	\qquad |t|<a.
	\]
	Here $ \mathcal{P}_{0\to t} $ denotes parallel transport along $\gamma$.
\end{proof}

\textbf{Acknowledgement.} I would thank Anton Petrunin for reading the first draft of this work.

\end{document}